\documentclass[11pt]{article}

\usepackage{amsmath}
\usepackage{amssymb}
\usepackage[notref,notcite]{}
\usepackage{hyperref}

\newenvironment{proofof}[1]{%
\noindent {\bf Proof of #1}}%
{\hspace*{\fill}$\Box$}

\newtheorem{theorem}{Theorem}[section]
\newtheorem{corollary} [theorem] {Corollary}
\newtheorem{definition} [theorem] {Definition}

\def\E{{\mathbb E}\,}
\def\Var{{\mathrm{Var}}\,}

\newcommand{\pr}{\mathbb P}

\title{Rearrangements of distributions on integers that minimize variance}
\author{Aistis Atminas 
\thanks{Xi'an Jiaotong-Liverpool University, 111 Ren'ai Road, Suzhou 215123, China. Email: Aistis.Atminas@xjtlu.edu.cn} \and Valentas Kurauskas
\thanks{Faculty of Mathematics and Informatics, Vilnius University, Naugarduko 24, LT-03225 Vilnius, Lithuania. Email: valentas.kurauskas@mif.vu.lt} }
\begin{document}

\maketitle

\begin{abstract}
Which permutations of a probability distribution on integers minimize variance?

Let $X$ be a random variable on a set of integers $\{x_1, \dots, x_N\}$
such that $\pr(X_i = x_i) = p_i$, $i \in \{1,\dots,N\}$. 
Let $(p^{(1)}, \dots, p^{(N)})$ be the sequence $(p_1, \dots, p_N)$
ordered non-increasingly. 
Let $X^+$ be the random variable
    defined by $\pr(X^+=0)=p^{(1)}$, $\pr(X^+=1) = p^{(2)}$, $\pr(X^+=-1)=p^{(3)}, \dots, \pr(X^+=(-1)^N \lfloor \frac {N} 2 \rfloor)=p^{(N)}$.
In this short note we generalize and prove the inequality $\Var X^+ \le \Var X$.
\end{abstract}

\section {Introduction}

Rearrangement inequalities, classically covered in Chapter X of Hardy, Littlewood and Pólya \cite{HLP}
have been applied to derive many other results, including isoperimetric inequalities, see, e.g., \cite{polyaszego}, and concentration function inequalities / variations of the Littlewood--Offord problem, see, e.g., \cite{tj2024, discreteinverse, lev1998, madimanwangwoo2018}.
Many of the latter results have a form similar to the following one.
Let $X_1, \dots, X_n$ be independent random variables supported on finite sets of integers, and let $X_1^+, \dots, X_n^+$ be
independent random variables with the corresponding rearranged distribution functions. Then there exist $a_1, \dots, a_n\in \{-1,1\}$ such that
\begin{equation}\label{eq.rearrangement}
    \max_{x \in \mathbb{Z}} \pr(X_1 + \dots + X_n = x) \le \max_{x \in \mathbb{Z}} \pr(a_1 X_1^+ + \dots + a_n X_n^+ = x).
\end{equation}
For example, Theorem~371 of~\cite{HLP} implies that (\ref{eq.rearrangement}) holds (with $a_1=1$, $a_2=-1$ and $a_3=a_4=\dots=1$) in the case when $X_i^+$ are symmetric for all $i \ge 3$  and the main result of \cite{lev1998} is that (\ref{eq.rearrangement}) holds when $X_i$ is distributed uniformly on a finite subset of $\mathbb{Z}$ (in this case the signs $a_i$ are not important).
 
Consider another particular case where $X_1, X_2, \dots$
are i.i.d. copies of an integer random variable $X$ with a finite support, and assume 
that the support of $X - k$ is not contained in~$s \mathbb{Z}$ for some integers $k$ and $s$, $s > 1$.
In this case the local limit theorem, see, e.g., Theorem 1 in Chapter VII of~\cite{petrov}, implies that
\[
    \max_{x \in \mathbb{Z}} \pr(X_1 + \dots + X_n = x) = \frac {1 + o(1)} { \sqrt{2 \pi n \Var X}}
\]
and so (\ref{eq.rearrangement}) holds for $n$ large enough (with $a_1=\dots=a_n=1$) if 
\begin{equation}\label{eq.var}
    \Var X^+ \le \Var X
\end{equation}
and if the equality in (\ref{eq.var}) is only achieved in the obvious cases when
$X - k \sim X^+$ or $X - k \sim - X^+$ for some integer $k$. 

The question whether (\ref{eq.var}) always holds arose while applying a similar argument in \cite{discreteinverse}.
In the present short note we provide a straightforward proof of (\ref{eq.var}) as we were not able to find it mentioned in the literature.
A special case of the present result is used, along with many other ideas, in the proof of a much more general result of \cite{discreteinverse}.

Let $f$ be the density function of an absolutely continuous random variable.
$f$ can be transformed, see Chapter 10.12 of \cite{HLP},
to obtain a 
density $f^*$
called the \emph{symmetric decreasing rearrangement} of $f$
which satisfies for any Borel set $B$ and the Lebesgue measure $\lambda$
\begin{equation}\label{eq.symdec}
\int_{[-\frac {|B|} 2, \frac {|B|} 2]} f^* d \lambda \ge \int_{B} f d \lambda.
\end{equation}
As for any non-negative random variable $\E X = \int_{t=0}^{\infty} \pr(X > t) dt$, for any $p \ge 1$ we 
have $\E |X - \E X|^p = \int_{t=0}^{\infty} p t^{p-1} \pr(|X - \E X| > t) dt$. If random variables $X$ and $X^*$ have densities $f$ and $f^*$ respectively, (\ref{eq.symdec}) implies
that $\pr(|X^*| > t) \le \pr(|X- \E X| > t)$ for any $t \ge 0$, so $\E |X^*|^p \le \E |X - \E X|^p$.
Thus, a `continuous' variant of (\ref{eq.var}),
as opposed to the integer variant that we consider here, follows rather easily,
and has been noted in the literature, see, e.g. \cite{wangmadiman2014}.


We will use the next definition.
\begin{definition}\label{def.D_f}
    Let $X$ be a random variable. Let $f: [0, +\infty) \to [0, +\infty)$ be a non-decreasing function. 
    Define a number
    \[
        D_f(X) := \inf_{a \in \mathbb{R}} \E f(|X - a|)
    \]
    and, in the case $D_f(X) < \infty$, the set
    \[
        M_f(X) := \{a \in \mathbb{R}: \E f(|X - a|) = D_f(X)\}.
    \]
\end{definition}

Thus every $f$ as above gives a measure of dispersion $D_f$
and a central tendency $M_f$. These statistics can also be generalized to the $d$-dimensional Euclidean space or other normed spaces.

\begin{theorem}\label{thm.main}
    Let $X$ be a random variable supported on a finite set of integers. Assume that $f: [0, +\infty) \to [0, +\infty)$
    is non-decreasing and $D_f(X) < \infty$. Then 
    \begin{equation}\label{eq.main}
        D_f(X^+) \le D_f(X).
    \end{equation}
    Furthermore, suppose that 
    $f$ is continuous
    with a
    positive derivative for $x > 0$ and a right derivative at $0$ such that  $f'(0+) = 0$.
    Then (\ref{eq.main}) is strict unless $X-k$ is distributed as $X^+$ or $-X^+$ for some integer $k$.
\end{theorem}
Recall that $m$ is a median of $X$ if $\pr(X \ge m) \ge \frac 1 2$ and $\pr(X \le m) \ge \frac 1 2$.
\begin{corollary}\label{cor.main}
    Let $X$ and $f$ be as in Theorem~\ref{thm.main}.
    \begin{enumerate}
        \renewcommand{\labelenumi}{\arabic{enumi})}
        \item If $f(x) = x$ for $x \ge 0$ then each value $m \in M_f(X)$ is a median of $X$
            and 
            $D_f(X) = {\rm MAD_{median}}(X) = \E |X - m|$, i.e., $D_f(X)$ is the mean absolute deviation
            of $X$ around the median.
        \item If $f(x) = x^2$ for $x \ge 0$ then $M_f(X) = \{\E X\}$ and $D_f(X) = \Var X$. 
    \end{enumerate}
    Thus ${\rm MAD_{median}}(X^+) \le {\rm MAD_{median}}(X)$ and $\Var X^+ \le \Var X$. Furthermore $\Var X^+ = \Var X$
    if and only if $X - k \sim X^+$ or $X-k \sim - X^+$ for some integer $k$.
\end{corollary}

\section {Proofs}

\begin{proofof}{Theorem~\ref{thm.main}}
Let $x_1, \dots, x_N$ be the support of $X$ listed
in such a way that the corresponding probabilities 
    $p^{(1)}, \dots, p^{(N)}$
    , given by $p^{(i)} = \pr(X = x_i)$,
are non-increasing.

It can be checked that 
$M_f(X)$ is non-empty when $f$ is continuous.
First assume that $M_f(X) \ne \emptyset$.
Let $a \in M_f (X)$. 
We denote $D_f(X) = {\bf p} \cdot {\bf v}$
where 
${\bf p} = (p^{(1)}, p^{(2)}, \ldots, p^{(N)})$ and
$${\bf v} = (f(|x_1-a|), f(|x_2-a|), \ldots, f(|x_N-a|)).$$

Let ${\bf{v'}}=(f(|x^{(1)}-a|), f(|x^{(2)}-a|), \ldots, f(|x^{(N)}-a|))$ be the sequence
$(f(|x_1-a|), f(|x_2-a|), \ldots, f(|x_N-a|))$ ordered non-decreasingly. 
Then, a classical result about the rearrangements of two sequences (e.g. Theorem~368 of~\cite{HLP}) implies that 
$${\bf p} \cdot {\bf v} \geq {\bf p} \cdot {\bf v'} . $$
Set $a' = \min(a-\lfloor a \rfloor, \lfloor a \rfloor +1 - a)$. In other words,
the number $a' \in [0, \frac{1}{2}]$ represents the distance between the number $a$ and its nearest integer.
Set $${\bf w}=(f(a'), f(1-a'), f(1+a'), f(2-a'), f(2+a'), \ldots, f(\lfloor\frac{N}{2}\rfloor+ (-1)^{N-1} a')).$$ 
Clearly, ${\bf w}$ is ordered non-decreasingly. Further, recalling that $\{x^{(1)}, \ldots, x^{(N)}\}$ is a set of $N$ distinct integers and $f$ is non-decreasing, it is not hard to see that every component of the vector $\bf{v'}-\bf{w}$ is 
non-negative.
Hence, we obtain that
 $${\bf p} \cdot {\bf v'} \geq {\bf p} \cdot {\bf w}.$$
Adding all the ingredients together we conclude that 
\begin{align}
D_f(X) &= {\bf p} \cdot {\bf v} \nonumber\\
&\geq {\bf p} \cdot {\bf v'} \nonumber\\
&\geq {\bf p} \cdot {\bf w}\label{eq.pw}\\
&= \E f(|X^+ - a'|) \label{eq.a'}\\
&\geq D_f(X^+)\nonumber .
\end{align}
This finishes the proof of (\ref{eq.main}) when $M_f(X) \ne \emptyset$.
For the general case, by definition, for any $\epsilon > 0$ we can find $a=a(\epsilon)$ such
that $\E f (|X-a|) \le D_f(X) + \epsilon$. The same argument as above shows that
$D_f(X) \ge D_f(X^+) - \epsilon$. Since $\epsilon > 0$ is arbitrary, 
we conclude that $D_f(X) \ge D_f(X^+)$.

Assume now the additional properties of $f$ stated in the second part of the theorem.
Now $f$ is continuous, so $M_f(X) \ne \emptyset$.
Assume $D_f(X) = D_f(X^+)$, but $X$ is not a translation of $X^+$ or $-X^+$. We will
follow the proof of (\ref{eq.main}) and obtain a contradiction.

Since translating by a constant does not change $D_f(X)$, we can assume without loss of generality that $\lfloor a \rfloor = 0$, equivalently, $a \in [0, 1)$.

When defining $\bf v$ and $x_1, \dots, x_N$ we may 
additionally assume
that $(\pr(X=x_1), -f(|x_1-a|))$, $\dots$, $(\pr(X=x_N), -f(|x_N-a|))$ is ordered non-increasingly in lexicographic order.

We claim that 
\begin{equation}\label{eq.vvw}
    {\bf v}={\bf v'}={\bf w}.
\end{equation}
To see the first equality, assume there exist $i$ and $j$
such that $i<j$ and $v_i > v_j$. Then due to to the ordering of $(x_i)$, it must be
$p^{(i)} > p^{(j)}$. This implies that
$p^{(i)} v_j + p^{(j)} v_i < p^{(i)} v_i + p^{(j)} v_j$,
so exchanging the atoms at $i$ and $j$ gives a random variable $X'$, with $D_f(X') < D_f(X)$,
which is a contradiction to  (\ref{eq.main}).

To see the second equality of (\ref{eq.vvw}), notice that since both of these vectors are ordered non-decreasingly,
if they are not equal, we must have that some component of $\bf{v'}-\bf{w}$ is positive, and hence (\ref{eq.pw}) is strict,
again 
a contradiction to  (\ref{eq.main}).

Suppose first that $a' \not \in \{0, \frac 1 2\}$.
Then, since $f$ is strictly increasing for $x > 0$, identity is the unique permutation that orders the components of $\bf v$ non-decreasingly.
When $a \in (0, \frac 1 2)$ this corresponds to placing the probabilities $p^{(1)}, \dots, p^{(N)}$
on $0, 1, -1, \dots$ respectively as in the distribution of $X^+$.
Similarly, when $a \in (\frac 1 2, 1)$, this corresponds to placing them on $1, 0, 2, -1, \dots$ respectively as in 
the distribution of $1 - X^+$. 

So we can assume that $a' \in \{0, \frac 1 2\}$. Then, if $a'=0$ we have $|x_{2k}|= |x_{2k+1}|$ for $k \in \{1,2, \dots\}$, and if $a'=\frac{1}{2}$ we have $|x_{2k-1} - a'|= |x_{2k} - a'|$ for $k \in \{1,2, \dots\}$.
It cannot be that for $a'=0$ we have
\begin{equation} \label{eq.symmetric2} 
 p^{(2k)} = p^{(2k+1)} \mbox{ for }k \in \{1,2,\dots\}
\end{equation}
or for $a'=\frac{1}{2}$ we have 
\begin{equation} \label{eq.symmetric}
      p^{(2k-1)} = p^{(2k)} \mbox{ for }k \in \{1,2,\dots\}
\end{equation}
since in these cases (\ref{eq.vvw}) implies that $X \sim X^+$ (the distribution is symmetric around $a'$). 

Suppose that $a' = a = 0$.
By the definition of $X^+$ we have $\pr(X^+=k) \ge \pr(X^+ = -k)$ for all $k \in \{1,2,\dots\}$.
Since (\ref{eq.symmetric2}) cannot hold, for some $k$ we have  $\pr(X^+=k) > \pr(X^+ = -k)$.
Consider the function $g(x) = \E f(| X^+ - x|)$.
By the assumptions on $f'$ of the theorem, we have
\begin{align*}
    & g'(0+) =  \pr(X^+=0) f'(0+)  -\sum_{k \in \mathbb{Z}\setminus\{0\}} {\rm sgn}(k) \pr(X^+ = k) f'(k) 
\\ &=  - \sum_{k \in \{1,2,\dots\}} (\pr(X^+ = k) - \pr(X^+ = -k)) f'(k) < 0,
\end{align*}
so 
$D_f(X^+) \le g(\delta) < g(0) \le D_f(X)$ for some $\delta > 0$, a contradiction.

Finally, suppose that $a' = a = \frac 1 2$.
Note that by the definition of $X^+$, $\pr(X^+ = 1-k) \ge \pr(X^+=k)$ for $k \in \{1,2,\dots\}$.
Since (\ref{eq.symmetric}) cannot hold, for some $k$ we have $\pr(X^+ = 1-k) > \pr(X^+ = k)$. Similarly as above 
\[
    g'\left(\frac 1 2\right) =  \sum_{k \in \{1,2,\dots\}} (\pr(X^+ = 1-k) - \pr(X^+ = k)) f'\left(k - \frac 1 2\right) > 0,
\]
so $D_f(X^+) \le g(\frac 1 2 - \delta) < g(\frac 1 2) = \E f(|X^+ - \frac 1 2|) \le D_f(X)$ for some $\delta > 0$,
again a contradiction.
\end{proofof}

\bigskip

\begin{proofof}{Corollary~\ref{cor.main}}
1) and 2) are folklore facts in statistics with straightforward proofs, see, e. g., Chapter 6 of~\cite{YK}. 
The conclusion follows by applying Theorem~\ref{thm.main}. Note that in 2) we have $f'(x) = 2 x > 0$ for $x > 0$
and $f'(0+)=0$ as required.
\end{proofof}

\end{document}